\documentclass[oneside,notitlepage,12pt]{article}

\usepackage{amssymb}
\usepackage[leqno]{amsmath}
\usepackage{amsfonts}
\usepackage{amsopn}
\usepackage{amstext}
\usepackage{amsthm}
\usepackage{enumitem}

\usepackage{verbatim}
\usepackage[colorlinks, backref]{hyperref}

\usepackage{fourier-orns}
\usepackage{clock} 

\providecommand{\cal}{\mathcal}
\renewcommand{\Bbb}{\mathbb}

\newenvironment{pf}{\begin{proof}}{\end{proof}}

\newcommand{\Aaa}{{\cal{A}}}
\newcommand{\Bee}{{\cal{B}}}
\newcommand{\Cee}{{\cal{C}}}
\newcommand{\Dee}{{\cal{D}}}
\newcommand{\Ef}{{\cal{F}}}

\newcommand{\Pee}{{\cal{P}}}

\newcommand{\Yu}{{\cal{U}}}

\newcommand{\Nat}{{\Bbb{N}}}

\newcommand{\sS}{{\mathcal{S}}}

\newcommand{\al}{\alpha}
\newcommand{\Gam}{\Gamma}

\newcommand{\sig}{\sigma}

\renewcommand{\phi}{\varphi}
\renewcommand{\rho}{\varrho}

\newcommand{\rest}{\restriction}

\newcommand{\ntr}{{n\in\omega}}

\newcommand{\loe}{\leq}
\newcommand{\goe}{\geq}

\newcommand{\subs}{\subseteq}
\newcommand{\sups}{\supseteq}
\newcommand{\nnempty}{\ne\emptyset}

\renewcommand{\iff}{\Longleftrightarrow}

\newcommand{\diam}{\operatorname{diam}}

\newcommand{\defi}{\,:=\,}

\newcommand{\Land}{\;\&\;}

\newtheorem{tw}{Theorem}[section]

\newtheorem{lm}[tw]{Lemma}
\newtheorem{prop}[tw]{Proposition}

\theoremstyle{definition}

\newtheorem{ex}[tw]{Example}

\newtheorem{pyt}[tw]{Question}

\theoremstyle{remark}

\newcommand{\setof}[2]{\{#1\mid #2\}}
\newcommand{\bigsetof}[2]{\Bigl\{#1\mid #2\Bigr\}}

\newcommand{\sett}[2]{\{#1\}_{#2}}
\newcommand{\sn}[1]{\{#1\}} 
\newcommand{\map}[3]{#1\colon #2 \to #3} 

\newcommand{\bF}{{\mathbb F}}

\providecommand{\nat}{\omega}

\newcommand{\separator}{\begin{center} \leafright \leafright \leafright \decotwo \decotwo \decotwo \leafleft \leafleft\leafleft
\end{center}}

\newcommand{\ci}{\operatorname{ci}}
\newcommand{\cir}{\operatorname{cir}}
\newcommand{\cu}{\operatorname{cu}}
\newcommand{\cur}{\operatorname{cur}}

\newcommand{\cB}{\mathcal B}
\newcommand{\cpl}{\mbox{\rm cpl}}
\newcommand{\cBc}{\mbox{\rm cpl}\mathcal B}
\newcommand{\bfind}[1]{\index{#1}{\em #1}}

\newcommand{\n}{\par\noindent}
\newcommand{\ssn}{\par\smallskip\noindent}
\newcommand{\mn}{\par\medskip\noindent}

\newcommand{\pars}{\par\smallskip}
\newcommand{\parm}{\par\medskip}

\title{Chain union closures}
\author{
{\sc Wies{\l}aw Kubi\'s}\footnote{The first author was supported by GA\v CR grant 20-22230L (Czech Science
Foundation).} \\
{\small Institute of Mathematics, Czech Academy of Sciences, Czech Republic}
\and
{\sc Franz-Viktor Kuhlmann}\footnote{The second author was supported by Opus grant 2017/25/B/ST1/01815 from the
National Science Centre of Poland.} \\
{\small Institute of Mathematics, University of Szczecin, Poland}
}

\date{April 1, 2024}

\begin{document}

\maketitle

\begin{abstract}
We study spherical completeness of ball spaces and its stability under expansions. We give some criteria
for ball spaces that guarantee that spherical completeness is preserved when the ball space is closed under
unions of chains. This applies in particular to the spaces of closed ultrametric balls in ultrametric spaces
with linearly ordered value sets, or more generally, with countable narrow value sets.
We show that in general, chain union closures of ultrametric spaces with partially ordered value
sets do not preserve spherical completeness.
{Further, we introduce and study the notions of chain union stability and of chain union rank, which measure how often the process of closing a ball space under all unions of chains has to be iterated until a ball space is obtained that is closed under unions of chains.}
	\ \\
	{\it MSC (2010):} 06A06, 54A05.

	\ \\
	{\it Keywords:} Ball space, ultrametric space, spherical completeness, narrow poset.
\end{abstract}

\tableofcontents

\section{Introduction}

In \cite{KK1,KK1a,KK2,KK3,KK4,KK5}, ball spaces are studied in order to provide a general framework for fixed
point theorems that in
some way or the other work with contractive functions. A {\em ball space} $(X,\cB)$ is a nonempty set $X$ together
with any nonempty collection of nonempty subsets of $X$. The completeness property necessary for the proof of
fixed point theorems is then encoded as follows. A \bfind{chain of balls} (also called a \emph{nest})
in $(X,\cB)$ is a nonempty subset of $\cB$ which is linearly ordered by inclusion. A ball space $(X,\cB)$ is
called \bfind{spherically complete} if every chain of balls has a nonempty intersection. Further, we
say that a ball space $(X, \cB)$ is \bfind{chain union closed} if the union of every chain in
$\cB$ is a member of $\cB$. We define $\cu(\cB)$ to be the family of all sets of the form $\bigcup \Cee$, where
$\Cee \subs \Bee$ is a chain (recall that, by default, chains of sets are supposed to be nonempty). More formally,
\[
\cu(\cB) = \bigsetof{ \bigcup \Cee }{\emptyset \ne \Cee \subs \cB, \text{ $\Cee$ is a chain} } \>.
\]
Hence a ball space $(X, \cB)$ is chain union closed if and only if $\cu(\cB) = \cB$. In the present paper, we
study the process of obtaining a chain union closed ball space from a given ball space and the question under
which conditions the spherical completeness of $(X,\cB)$ implies the spherical completeness of $(X,\cu(\cB))$.

A ball space $(X, \cB)$ is said to be \bfind{chain union stable} if $(X, \cu(\cB))$ is chain union closed. In other words, for every chain $\Dee \subs \cu(\cB)$ there exists a chain $\Cee \subs \cB$ with $\bigcup \Cee = \bigcup \Dee$.
{Clearly, every chain union closed space is chain union stable. Furthermore, $(X, \cB)$ is chain union stable if and only if $\cu(\cB)$ is chain union closed.}

\begin{ex}
	Let $\cB$ be the family of all finite nonempty subsets of a fixed set $X \nnempty$. Then $\cu(\cB)$ is the family of all nonempty countable subsets of $X$. Note that $(X, \cB)$ is chain union stable if and only if $X$ is countable.
\end{ex}


\separator

The main inspiration for these definitions and questions is taken from the theory of ultrametric spaces and their
ultrametric balls. An \bfind{ultrametric} $u$ on a set $X$ is a function from $X\times X$ to a partially ordered
set $\Gamma$ with smallest element $\bot$, such that for all $x,y,z\in X$ and all $\gamma\in\Gamma$,
\begin{enumerate}
	\item[(U1)] $u(x,y)=\bot$ if and only if $x=y$,
	\item[(U2)] if $u(x,y)\leq\gamma$ and $u(y,z)\leq\gamma$, then
	$u(x,z)\leq\gamma$,
	\item[(U3)] $u(x,y)=u(y,x)$ \ \ \ (symmetry).
\end{enumerate}
Condition (U2) is the ultrametric triangle law; if $\Gamma$ is linearly ordered, it
can be replaced by
\begin{enumerate}
	\item[(UT)] $u(x,z)\leq\max\{u(x,y),u(y,z)\}$.
\end{enumerate}
{When dealing with such ultrametric spaces, we can say that a set $A$ has \emph{diameter} $\loe \gamma$ if $u(x,y) \loe \gamma$ for every $x,y \in A$. On the other hand, the diameter of $A$ may not be defined, unless the value set $\Gamma$ is a complete meet semilattice.}
A \bfind{closed ultrametric ball} is a set $B_\alpha(x) \defi \setof{y\in X}{u(x,y)\leq \alpha}$, where $x\in X$
and $\alpha\in\Gamma$. The problem with general ultrametric spaces is that closed balls $B_\alpha(x)$ are not
necessarily precise, that is, there may not be any $y\in X$ such that $u(x,y)=\alpha$. Therefore, we prefer to
work only with \bfind{precise} ultrametric balls, which we can write in the form
$$
B(x,y)\defi \setof{z\in X}{u(x,z)\leq u(x,y)},
$$
where $x,y\in X$.
{Note that a precise ultrametric ball $B(x,y)$ has diameter precisely $u(x,y)$.}
We obtain the \bfind{ultrametric ball space} $(X,\cB_u)$ from $(X,u)$
by taking $\cB_u$ to be the set of all such balls $B(x,y)$.
Specifically, $\cB_u \defi \setof{B(x,y)}{x,y \in X}$.

More generally, an ultrametric ball is a set
\[
B_S (x) \>:=\> \setof{y\in X}{u(x,y)\in S}\>,
\]
where $x\in X$ and $S$ is
an initial segment of $\Gamma$. We call $X$ together with the collection of all ultrametric balls
the \bfind{full ultrametric ball space} of $(X,u)$. Every ultrametric ball can be written as the union over a chain of precise balls:
\[
B_S (x) \>=\> \bigcup\{B(x,y)\mid u(x,y)\in S\}\>.
\]
Hence the full ultrametric ball space is just $(X,\cu(\cB_u))$.

Typically, ultrametric spaces are considered with a linearly ordered value set, in which case the ball structure is a tree, in the sense that given a ball $B$, the set $\setof{C \sups B}{C \text{ is a ball }}$ is a chain. This is not true when the distance set is partially ordered, however we still have the following easy and well-known weaker fact.

\begin{prop}\label{PROPkjsdvah}
	Let $(X,u)$ be an ultrametric space and let $B_0$, $B_1$ be ultrametric balls with nonempty intersection and the same diameter. Then $B_0 = B_1$.
\end{prop}

We will exploit the fact that classical ultrametric spaces are tree-like.
A ball space $(X,\cB)$ is called \bfind{tree-like} if for
every $B_1, B_2 \in \cB$ the following implication holds.
\begin{equation}
	B_1 \cap B_2 \nnempty \implies B_1 \subs B_2 \text{ or } B_2 \subs B_1.
	\tag{I}\label{eqaI}
\end{equation}
See \cite{KKub} for some remarks on tree-like ball spaces.

\separator

We now formulate our main results.

\begin{tw}                                  \label{MT1}
	Let $(X,\cB_u)$ be the ball space of an ultrametric space $(X,u)$ with linearly ordered value set. Then the following
	assertions hold:
	\begin{enumerate}[itemsep=0pt]
		\item[{\rm(1)}] The ball space $(X,\cB_u)$ is chain union stable.
		\item[{\rm(2)}] If $(X,\cB_u)$ is spherically complete, then so is $(X,\cu(\cB_u))$.
	\end{enumerate}
\end{tw}

This result is a consequence of Theorem~\ref{tlbs} and the fact that ultrametric spaces with linearly ordered value set are tree-like.

The next result is actually a simple example, therefore we present the proof immediately.
Given two ball spaces $(X,\cB)$ and $(X,\cB')$ on the same set $X$, we call $(X,\cB')$ an \bfind{expansion}
of $(X,\cB)$ if $\cB\subseteq \cB'$. In general, we cannot expect the existence of chain
union closed expansions which preserve spherical completeness:

\begin{tw}\label{MT3}
	There exists a countable spherically complete ultrametric space with a countable partially ordered value set, whose
	ultrametric ball space does not admit any expansion that is chain union closed and spherically complete.
\end{tw}

\begin{pf}
	Let $\Gam = \sn \emptyset \cup 
	\Gam'$, where $\Gam'$ is the set of all closed intervals $[k,\ell]$, where $k < \ell$ are nonnegative integers. Then $(\Gam, \subs)$ is a countable poset and $\emptyset$ is its smallest element. Define $\map u {\Nat^2}\Gam$ by $u(x,y) = u(y,x) = [x,y]$ if $x < y$
	and $u(x,x) = \emptyset$. Clearly, $(\Nat, u)$ is a spherically complete ultrametric space; its balls are finite intervals of natural numbers. Now suppose $\cB'$ is a chain union closed expansion of $\cB_u$.
	Then for each $n \in \Nat$ the unbounded interval
	$$[n,\infty)\cap \Nat = \bigcup_{m>n}[n,m]\cap \Nat$$
	belongs to $\cB'$. Finally, $\sett{[n,\infty)}{n \in \Nat}$ is a chain in $\cB'$ whose intersection is the empty set.
\end{pf}

So we cannot hope for extending Theorem~\ref{MT1}(2) to partially ordered value sets. It turns out that part (1) of this theorem holds in a more general setting.

\begin{tw}\label{THMJedenTrzi}
	Assume $(X,u)$ is an ultrametric space with a countable partially ordered value set. Then its ball space $(X, \cB_u)$ is chain union stable. 
\end{tw}

We shall deduce this result from a more abstract statement involving the concept of chain regularity, inspired by a similar notion in category theory.

We finish this section with the following simple example exhibiting problems with uncountable cardinalities.

\begin{ex}
	Let $X = \omega \times \omega_1$ and let $\cB$ consist of all rectangles $[0,n) \times [0,\alpha)$ where both $n \in \omega$, $\al \in \omega_1$ are positive. Note that $\cu(\cB)$ consists of all rectangles of the form $[0,\xi) \times [0,\beta)$, where either $\xi<\omega$, $\beta<\omega_1$ or $\xi = \omega$ and $\beta < \omega_1$ or $\xi < \omega$ and $\beta = \omega_1$.
	In particular, $X \notin \cu(\cB)$.
	It follows that $(X, \cB)$ is not chain union stable.
\end{ex}

\section{Chain regularity}

We say a ball space $(X, \cB)$ is \bfind{chain regular} if for every chain $\Cee \subs \cB$, for every ball $B \subs \bigcup \Cee$ there exists $C \in \Cee$ with $B \subs C$.
This notion is well justified by the following:

\begin{prop}\label{PROPbvksdbv}
	Every ultrametric ball space (with a partially ordered value set) is chain regular.
\end{prop}

\begin{pf}
	Assume $A \subs \bigcup \Cee$, where $\Cee$ is a chain of precise ultrametric balls. Let $a, b \in A$ be such that $A = B(a,b)$ and choose $C \in \Cee$ such that $a,b \in C$. Then the diameter of $C$ is $\goe$ the diameter of $A$, therefore $A \subs C$.
\end{pf}

\begin{lm}\label{LMfisgbuisgh}
	Assume $(X,\cB)$ to be a chain regular ball space such that all chains in $\cu(\cB)$ have countable cofinality. Then $(X,\cB)$ is chain union stable.
\end{lm}

\begin{pf}
	Let $\Cee = \sett{C_n}{\ntr}$ be a chain in $\cu(\cB)$. For each $\ntr$ choose a chain $\Aaa_n \subs \cB$ whose union is $C_n$. By our assumptions, $\Aaa_n$ has countable cofinality.

	Using standard induction together with chain regularity, construct a matrix of balls $\sett{B_{i,j}}{i,j < \nat}$ such that
	\begin{enumerate}[itemsep=0pt]
		\item $B_{i,j} \in \Aaa_i$ and $\bigcup_{\ntr} B_{i,n} = C_i$,
		\item $B_{i,j} \subs B_{i+1,j}$,
	\end{enumerate}
	for every $i,j < \nat$.
	Finally, the union of the diagonal $\sett{B_{n,n}}{\ntr}$ is $\bigcup_{\ntr} C_n$.
\end{pf}

We shall need the following simple result concerning partially ordered sets.

\begin{lm}\label{LMerighworug}
	Let $(P ,\loe)$ be a partially ordered set, $Q \subs P$ and $C$ a chain in $P$ such that each element of $C$ is the supremum of a chain in $Q$. Let $\kappa$ be an uncountable regular cardinal such that $|Q| < \kappa$. Then $C$ does not contain a copy of $(\kappa, \in)$.
\end{lm}

\begin{pf}
	Assume $\map \phi {\kappa}C$ preserves $\loe$. We are going to show that $\phi$ is eventually constant.
	Refining\footnote{By \emph{refining} a set we mean removing ``unnecessary'' or ``irrelevant'' elements.} $Q$, we may assume that for each $q\in Q$ there is $\al_q < \kappa$ with $q  \loe \phi(\al_q)$.
	Let $\delta = \sup \setof{\al_q}{q \in Q}$. Then $\delta<\kappa$, because $\kappa$ is regular. Furthermore $\phi(\delta)$ is an upper bound of the whole set $Q$, therefore $\phi(\beta) = \phi(\delta)$ for every $\beta>\delta$.
\end{pf}

\begin{pf}[Proof of Theorem~\ref{THMJedenTrzi}]
	We can apply Lemma~\ref{LMfisgbuisgh}, as long as we prove that all chains in $\cu(\cB_u)$ have countable cofinality.

	Assume $\Cee = \sett{C_\al}{\al<\omega_1}$ is an $\omega_1$-chain in $\cu(\cB_u)$, i.e., $C_\al \subs C_\beta$ whenever $\al<\beta$.
	Choose $z \in C_0$. 
	Each $C_\al$ is the union of a chain of balls $\sett{B_{\al,n}}{\ntr} \subs \cB_u$ and we may assume that $z \in B_{\al,0}$ for every $\al <\omega_1$. We claim that the family $\Ef = \setof{B_{\al,n}}{\al< \omega_1, \ntr}$ is countable. Then Lemma~\ref{LMerighworug} gives $C_\beta = C_\delta$ for all big enough $\beta, \delta < \omega_1$.
	
	Given $\al, \al'<\omega_1$, $n,n'<\omega$, note that $B_{\al,n} = B_{\al', n'}$ whenever $\diam(B_{\al,n}) = \diam(B_{\al',n'})$ (Proposition~\ref{PROPkjsdvah}). This shows that $\Ef$ is countable, which completes the proof.
\end{pf}


\separator

The structure of ultrametric spaces with partially ordered value sets is in general much more complex than in the
case of linearly ordered value sets. But we can at least prove the following.
Recall that a partially ordered set (``poset'') is \bfind{narrow} if it contains no infinite sets of pairwise
incomparable elements.

\begin{tw}\label{MT2}
Let $(X,\cB_u)$ be the ball space of an ultrametric space with countable narrow value set. Then $(X,\cB_u)$ is chain union stable.
\end{tw}

\noindent
This theorem will be deduced from Theorem~\ref{ThmHevyNarow} at the end of Section~\ref{Sectcucb}, where we
study chain union closures for a class of ball spaces that constitutes a slight generalization of
the class of ball spaces whose chain intersection closures we have studied in \cite[Section~3]{KKub}.

We do not know whether the countability assumption can be dropped in Theorem~\ref{MT2}.

\section{Chain union and chain intersection closures}             \label{Sectcuc}
Let $\cB$ be a nonempty family of nonempty sets. Using transfinite recursion, we define $\cu_\al(\cB)$ and
$\ci_\al(\cB)$ for each ordinal $\al$, as follows:
\begin{eqnarray*}
\cu_0(\cB) &=& \cB, \qquad \cu_\al(\cB) \>=\> \cu\left( \bigcup_{\xi < \al} \cu_\xi(\cB) \right) \text{ for }
\al > 0\>,\\
\ci_0(\cB) &=& \cB, \qquad \ci_\al(\cB) = \ci\left( \bigcup_{\xi < \al} \ci_\xi(\cB) \right) \text{ for }
\al > 0\>,
\end{eqnarray*}
where $\ci(\Bee)$ is defined in~\cite{KKub} as 
\[
\ci(\cB) = \bigsetof{ \bigcap \Cee }{\emptyset \ne \Cee \subs \cB, \text{ $\Cee$ is a chain} } \>.
\]
We have $\cu(\cB)=\cu_1(\cB)$ and $\ci(\cB)=\ci_1(\cB)$.
We observe that
\begin{equation}                   \label{cicusubs}
\cB\subs\cB'\>\Rightarrow\> \ci_\al(\cB)\subs\ci_\al(\cB')\mbox{ and }\cu_\al(\cB)\subs\cu_\al(\cB')
\mbox{ for all }\al\>.
\end{equation}

We define the \bfind{chain union rank} of $\cB$, denoted by $\cur(\cB)$, to be the smallest ordinal $\al$ such
that $\cu_{\al+1}(\cB) = \cu_\al(\cB)$.
Thus, $\cur(\cB) = 0$ if and only if $\cB$ is chain union closed, while $\cur(\cB) \loe 1$ means that in order to
make $\cB$ chain union closed, it suffices to extend it by adding all unions of chains.
In general, we call $(X, \cu_\al(\cB))$, with $\al = \cur(\cB)$, the \bfind{chain union closure} of $(X, \cB)$. It
could also be described as a ball space $(X, \cB')$, where $\cB' \sups \cB$ is minimal such that $\cB'$ is stable
under unions of chains.

Recall from the introduction that a ball space $(X,\cB)$ is \bfind{chain union stable} if for every nonempty family $\bF$ consisting of
chains in $\cB$ such that $\setof{\bigcup \Cee}{\Cee \in \bF}$ is a chain, there exists a chain $\Yu \subs \cB$
satisfying
\[
\bigcup \Yu \>=\> \bigcup_{\Cee \in \bF} \bigcup \Cee\>.
\]
Obviously, this is equivalent to saying that $\cur(X, \Bee) \loe 1$.
In \cite{KKub}, the notions of \bfind{chain intersection rank}, denoted by $\cir(\cB)$, of \bfind{chain
intersection closure} and of \bfind{chain intersection stable} were defined analogously, just with ``cu''
replaced by ``ci''.
Note that both $\cur$ and $\cir$ are well-defined, because there is a cardinality on the powerset of the ball space.

\parm
If $(X,\cB)$ is a ball space, then also $(X,\cBc)$ with
\[
\cBc\>:=\>\{X\setminus B\mid B\in\cB, B\ne X\}
\]
is a ball space; we will call it the \bfind{complement ball space} of $(X,\cB)$. Note that $\cBc=\mbox{cpl}
(\cB\setminus\{X\})$. We collect a number of useful properties of chain unions and chain intersections.
\begin{lm}
Take a ball space $(X,\cB)$.\ssn
1) If $\sS\subs \cB$ is a finite set and $\cB\setminus \sS\ne\emptyset$, then
\ssn
i) $\ci_\al(\cB) \>=\>\ci_\al(\cB\setminus \sS)\cup \sS$ and
$\,\cu_\al(\cB) \>=\>\cu_\al(\cB\setminus \sS)\cup \sS$ for all $\al$,\n
ii) $\cir(\cB)\leq\cir(\cB\setminus \sS)$ and $\,\cur(\cB)\leq\cur(\cB\setminus \sS)$.
\mn
2) If $X\in\cB$ and $\cB\setminus\{X\}\ne\emptyset$, then $\ci_\al(\cB)\setminus\ci_\al(\cB\setminus \{X\})=
\{X\}$ for all $\al$ and $\cir(\cB)=\cir(\cB\setminus \{X\})$.
\mn
3) For all $\al$, $(X,\ci_\al(\cBc))$ is the complement ball space of $(X,\cu_\al(\cB))$.
\mn
4) We have that
\[
\cur(\cB)\>\geq\>\cir(\cBc)\>.
\]
If $\cir(\cBc)=\beta$ and $X\in\cu_\beta(\cB)$, then $\cur(\cB)=\cir(\cBc)$.
\end{lm}
\begin{pf}
1): We treat the case of chain intersections; the case of chain unions is analogous. We have that
$\ci(\cB)=\ci(\cB\setminus \sS)\cup \sS$ as all members of $\sS$ can be removed from any infinite chain without
changing the intersection. This implies assertion i) for $\al=1$ in the case of chain intersections.

Now we proceed by induction on $\al$. Assume that assertion i) holds for $\ci$ and $\ci_\al\,$. Then
\begin{eqnarray*}
\ci_{\al+1}(\cB) &=& \ci(\ci_\al(\cB)) \>=\> \ci(\ci_\al(\cB\setminus \sS)\cup \sS)\>=\>
\ci((\ci_\al(\cB\setminus \sS)\cup \sS)\setminus \sS)\cup \sS \\
&=& \ci(\ci_\al(\cB\setminus \sS))\cup \sS\>=\>\ci_{\al+1}(\cB\setminus \sS)\cup \sS\>,
\end{eqnarray*}
where we have used our assertion for $\ci_\al$ for the second equality,
our assertion for $\ci$ for the third equality, {and then again our assertion for $\ci$ with $\ci_\alpha(\cB\setminus \sS)$ in place of $\cB$ for the fourth equality.}
This proves the successor case of the induction. The limit case is straightforward.

\pars
In order to prove assertion ii), assume that $\cir(\cB\setminus \sS)=\al$, that is,
$\ci_{\al+1}(\cB\setminus \sS))=\ci_\al(\cB\setminus \sS)$. Then by assertion i),
\[
\ci_{\al+1}(\cB)\>=\>\ci_{\al+1}(\cB\setminus \sS)\cup \sS\>=\>\ci_{\al}(\cB\setminus \sS)\cup \sS
\>=\>\ci_{\al}(\cB)\>.
\]
This proves that $\cir(\cB)\leq\al=\cir(\cB\setminus \sS)$.

\mn
2): Since the only chain that has $X$ as its intersection is $\{X\}$, no ball space $\cB$ satisfies
$X\in\ci(\cB\setminus\{X\})$. By induction, $X\notin\ci_\al(\cB\setminus\{X\})$ for all $\al$. Hence if
$X\in\cB$ and $\cB\setminus\{X\}\ne\emptyset$, then $X\in\ci_\al(\cB)\setminus\ci_\al(\cB\setminus \{X\})$, and it
follows from assertion i) of part 1) that $\ci_\al(\cB)\setminus\ci_\al(\cB\setminus \{X\})=\{X\}$.

From assertion ii) of part 1) we know that $\cir(\cB)\leq\cir(\cB\setminus \{X\})$; we have to show that
also ``$\geq$'' holds. Assume that $\ci_{\al+1}(\cB)=\ci_\al(\cB)$. Then by what we have just proved,
\[
\ci_{\al+1}(\cB\setminus \{X\})\>=\>\ci_{\al+1}(\cB)\setminus \{X\}\>=\>\ci_{\al}(\cB)\setminus \{X\}
\>=\>\ci_{\al}(\cB\setminus \{X\})\>.
\]
This proves the desired inequality and thus the second assertion of part 2).

\mn
3): The assertion is proven by induction on $\al$ using the fact that the complement of the union of a chain
$\{B_i\}_{i\in I}$ is the intersection of the chain $\{X\setminus B_i\}_{i\in I}$.

\mn
4): Let $\cur(\cB)=\al$, that is, $\cu(\cu_\al(\cB))=\cu_\al(\cB)$. Pick a chain $\Cee$ in $\ci_\al(\cBc)$
such that $\bigcap\Cee\ne\emptyset$. By part 3), $\{X\setminus B\mid B\in \Cee\}$ is a subset of $\cu_\al(\cB)$,
and it is also a chain. By assumption, $B':=\bigcup\{X\setminus B\mid B\in \Cee\}\in \cu_\al(\cB)$. Since
$\bigcap\Cee\ne\emptyset$, we have that $B'\ne X$. Using part 3) again, $\bigcap\Cee=X\setminus B'\in
\ci_\al(\cBc)$. We have proved that $\ci_\al(\cBc)$ is chain intersection closed, which shows that $\cir(\cBc) \loe \al$. Hence our first assertion holds.

Now assume that $\cir(\cBc)=\beta$ and that $X\in\cu_\beta(\cB)$.  By what we have proved before, it suffices to show
that $\cur(\cB)\leq \beta$. Pick a chain $\Cee$ in $\cu_\beta(\cB)$; we wish to show that $\bigcup\Cee\in
\cu_\beta(\cB)$. As $X\in\cu_\beta(\cB)$, we may assume that $\bigcup\Cee\ne X$, so that $B':=\bigcap\{X\setminus B
\mid B\in \Cee\}\ne\emptyset$. By part 3), $\{X\setminus B\mid B\in \Cee\}$ is a subset of $\ci_\beta(\cBc)$,
and it is also a chain. Since $\cir(\cBc)=\beta$, we find that $B'\in \ci_\beta(\cBc)$. Using part 3) again,
$\bigcup\Cee=X\setminus B'\in\cu_\beta(\cB)$. This shows that $\cur(\cB)\leq \beta = \cir(\cBc)$, as desired.
%
\end{pf}
Note that it can happen that
\[
\cir(\cBc)\>=\>\cur(\cB)\><\>\cur(\cB\setminus \{X\})\>.
\]
For example,
take $X=\Bbb N$ and $\cB$ to be the collection of all initial segments of $\Bbb N$. Then $(X,\cB)$ is chain union
closed, while $\cur(\cB\setminus \{X\})=1$. Further, we see that $\cpl(\cB)$ is the collection of
all final segments of $\Bbb N$. It is chain intersection closed, i.e., $\cir(\cpl(\cB))= \cur(\cB)=0$.

\par\medskip
We will now demonstrate by an example that both the chain union rank and the chain intersection rank of a ball
space can
be equal to any ordinal $\al$. Since the ball space $(X,\Pee(X)\setminus\{\emptyset\})$ for nonempty $X$ is both
chain union and chain intersection closed, we have to show this only for the case of $\al\geq 1$.
\begin{ex}\label{ExNoWERTwegho}
Take an ordinal $\al\geq 1$ and set $X:=\aleph_\al$.
For $\beta$ any ordinal, define $\cB_\beta$ to be the collection of all nonempty subsets of $\aleph_\al$ of
cardinality smaller than or equal to $\aleph_\beta\,$. Set $\cB:=\cB_0\,$. We note that $X\notin\cB$ since
$\al\geq 1$.
A standard transfinite induction shows that $\cu_\beta(\cB) = \cB_\beta$ for every $\beta \loe \al$.


\parm
Since the subsets of $\aleph_\al$ have cardinality at most $\aleph_\al\,$, it follows that
$\cu_\al(\cB)=\cB_\al=\Pee(\aleph_\al)\setminus\{\emptyset\}=\Pee(X)\setminus\{\emptyset\}$. Therefore,
$\cu_\al(\cB)$ is chain union closed. On the other hand,
$\cu_\beta(\cB)$ is not chain union closed for any $\beta<\al$. Hence, $\cur(\cB)=\al$.

\pars
Finally, we show that $\cir(\cBc)=\al$. By part 4) of the preceding lemma, $\beta:=\cir(\cBc)\leq\cur(\cB)=\al$.
Applying part 1)i) of the same lemma with $S=\{X\}$, we obtain that $\cu_\beta(\cB\cup\{X\})=
\cu_\beta(\cB)\cup\{X\}$. We have that
\[
\beta\>=\>\cir(\cBc)\>=\>\cir(\cpl(\cB\cup\{X\}))\>=\>\cur(\cB\cup\{X\})\>,
\]
where the last equality follows from part 4) of the previous lemma. This implies that $\cu_\beta(\cB)\cup\{X\}
=\cB_\beta\cup\{X\}$ is chain union closed. But as we have seen above, this can only be if $\cB_\beta\cup\{X\}=
\cB_\al\,$. Since $X$ is not the only subset of $X$ of cardinality $\aleph_\al$, this can only be the case if
$\beta=\alpha$, showing that $\cir(\cBc)=\al$.
\end{ex}

Modifying the example above, namely, declaring $\cB$ to be the family of all non\-empty finite sets, we obtain a spherically complete ball space $(X,\cB)$ such that $(X, \cu(\cB))$ is not spherically complete.

\section{Chain union closed or stable ball spaces}                      \label{Sectcucb}
We wish to find conditions for a ball space $(X,\cB)$ to be chain union closed/stable.

\begin{tw}\label{tlbs}
Let $(X, \cB)$ be a tree-like ball space. Then $\cu(\cB)$ is tree-like and chain union closed. Furthermore, if $(X, \cB)$ is spherically complete then so is $(X, \cu(\cB))$.
\end{tw}

\begin{pf}
	Fix $D_0, D_1 \in \cu(\cB)$ with $a \in D_0 \cap D_1$. Let $\Cee_0$, $\Cee_1$ be chains in $\cB$ with $D_i = \bigcup \Cee_i$, $i<2$.
	We may assume $a \in \bigcap \Cee_i$, by refining the chain. Now $\Cee_0 \cup \Cee_1$ is a chain, because $\cB$ is tree-like.
	Suppose $D_1 \not\subs D_0$ and fix $b \in D_1 \setminus D_0$. Choose $C_1 \in \Cee_1$ with $b \in C_1$. Now $C \subs C_1$ for every $C \in \Cee_0$, because the other inclusion is impossible. Hence $D_0 = \bigcup \Cee_0 \subs C_1 \subs \bigcup \Cee_1 = D_1$. This shows that $\cu(\cB)$ is tree-like.
	
	Now, take a chain $\Dee$ in $\cu(\cB)$ and a point $a \in \bigcup \Dee$. After refining $\Dee$, we may assume $a \in D$ for all $D \in \Dee$. Next, for $D \in \Dee$, let $D = \bigcup \Cee_D$, where $\Cee_D$ is a chain in $\cB$. Again, we may assume that $a \in C$ for every $C \in \Cee_D$. Since $(X, \cB)$ is tree-like, we conclude that $\bigcup_{D \in \Dee} \Cee_D$ is a chain in $\cB$ with the same union as $\Dee$. Hence, $\cu(\cB)$ is chain union closed.
	
	Finally, assume $(X, \cB)$ to be spherically complete and fix a strictly decreasing chain $\sett{D_\al}{\al<\kappa}$ in $\cu(\cB)$.
	Fix $\al<\kappa$. Choose $a \in D_{\al+1}$ and $b \in D_\al \setminus D_{\al+1}$. Since $D_\al$ is the union of a chain from $\cB$, we can find $B_\al$ in that chain such that $a, b \in B_\al$. Hence $B_\al \subs D_{\al}$. We also have $D_{\al+1} \subs B_\al$, because $a \in D_{\al+1} \cap B_\al$ and $b \in B_\al \setminus D_{\al+1}$, so the opposite inclusion is impossible. Finally, $\sett{B_\al}{\al<\kappa}$ is a chain in $\cB$ with $\bigcap_{\al<\kappa} B_\al = \bigcap_{\al<\kappa} D_\al$.
\end{pf}

Every ultrametric space with linearly ordered value set is tree-like, since in this case property (\ref{eqaI})
follows from the ultrametric triangle law. Hence Theorem~\ref{MT1} follows from Theorem~\ref{tlbs}.

\separator

For our next theorem, we will need two lemmas that reflect important and well known properties of narrow posets. The first one is an immediate consequence of the Erd\H{o}s-Dushnik-Miller Theorem.
\begin{lm}[{cf.~\cite[Lemma 3.2]{KKub}}]\label{LmDusek}
Let $(\Pee,\loe)$ be a narrow poset and  $\Aaa \subs \Pee$ infinite. Then there exists a chain $\Cee \subs \Aaa$
such that $|\Cee| = |\Aaa|$.
\end{lm}

Recall that a subset $A$ of a poset $(P,\loe)$ is \bfind{directed} if for every $a_0,a_1 \in A$ there is $b \in A$
with $a_0 \loe b$ and $a_1 \loe b$. The following fact can be found in~\cite{Fraisse}. We presented a proof in
\cite[Lemma~3.3]{KKub}.
\begin{lm}\label{LmFinDirek}
Every narrow poset is a finite union of directed subsets.
\end{lm}

We need one more general combinatorial fact.

\begin{lm}\label{LMbergibrgh}
	Assume $\Aaa = \Aaa_0 \cup \dots \cup \Aaa_{k-1}$ is a family of sets such that $\bigcup \Aaa$ is the union of a chain of sets that are unions of chains in $\Aaa$. Then there exists $j < k$ such that $\bigcup \Aaa = \bigcup \Aaa_j$.
\end{lm}

\begin{pf}
	Our assumption says that $\bigcup \Aaa = \bigcup_{\al<\kappa}\Cee_\al$ where $\Cee_\al \subs \Aaa$ is a chain for each $\al<\kappa$ and $\bigcup \Cee_\al \subs \bigcup \Cee_\beta$ whenever $\al<\beta$. We may assume that $\kappa$ is a regular cardinal.
	For each $\al<\kappa$ we choose $j(\al)<k$ such that $\bigcup \Cee_\al = \bigcup (\Cee_\al \cap \Aaa_{j(\al)})$. If $\kappa$ is finite then $\bigcup \Aaa = \bigcup(\Cee_{(\kappa-1)} \cap \Aaa_{j(\kappa-1)}) \subs \bigcup \Aaa_{j(\kappa-1)}$.
	Otherwise, there is an unbounded set $S \subs \kappa$ such that $j(\al) = j$ for every $\al \in S$.
	Then 
	$$\bigcup \Aaa = \bigcup_{\al \in S} \Cee_\al \subs \bigcup \Aaa_j.$$
\end{pf}

\par\medskip
We will now extend Theorem~\ref{tlbs} to a larger class of ball spaces $(X,\cB)$. For every $z\in X$,
we set
\[
\cB(z)\>:=\>\{B\in\cB\mid z\in B\}\>.
\]
\begin{tw}                                          \label{ThmHevyNarow}
Let $(X,\cB)$ be a ball space such that for every $z\in X$, the poset $(\cB(z),\subs)$ is narrow and admits only
countable strictly increasing sequences. Then $(X,\cB)$ is chain union stable.
\end{tw}

\begin{pf}
	Fix a strictly increasing chain $\Dee = \sett{D_\al}{\al<\kappa}$, where $\kappa$ is an infinite regular cardinal. We may assume that for each $\al<\kappa$ there is some $z_\al \in D_\al \setminus \bigcup_{\xi<\al} D_\xi$ (otherwise, replace $D_\al$ by $D_{\al+1}$). Let $z := z_0$.
	Each $D_\al$ is the union of a chain from $\cB(z)$, so we may restrict attention to the narrow poset $(\cB(z), \subseteq)$ in which all chains have countable cofinality.
	Using transfinite induction, we construct for each $\alpha<\kappa$ a chain $\sett{B_{\alpha,n}}{\ntr} \subs \cB(z)$ with $\bigcup_{\ntr}B_{\alpha,n} = D_\alpha$ and $z_\al \in B_{\alpha,0}$.
	By this way, if $B_{\al,n} \subs B_{\al',n'}$ then $\al \loe \al'$, because otherwise $z_\al \in B_{\al,n} \setminus B_{\al',n'}$.
	It follows that $\Aaa := \sett{B_{\al,n}}{\al<\kappa, n<\omega}$ is well-founded.
	
	By Lemma~\ref{LmDusek}, we conclude that $\Aaa$ is countable and consequently $\kappa = \aleph_0$.
	Finally, Lemma~\ref{LmFinDirek} says that $\Aaa$ is the union of finitely many directed subfamilies. The union of one of them must be equal to $\bigcup \Dee$, according to Lemma~\ref{LMbergibrgh}. A countable directed family of sets obviously has a cofinal chain, therefore $\bigcup \Dee$ is the union of an $\omega$-chain of balls.
\end{pf}

\noindent
We need the following lemma in order to deduce Theorem~\ref{MT2}.

\begin{lm}\label{LMetighueri}
Let $(X,\cB)$ be the ball space of an ultrametric space with countable narrow value set $\Gamma$. Then for each
$z\in X$, $(\cB(z),\subs)$ is narrow and admits only countable strictly increasing sequences.
\end{lm}
\begin{pf}
Take $\{B(x_i,y_i)\}_{i<\omega} \subs \cB(z)\,$. Then there are $k < \ell < \omega$ such that $u(x_k,y_k) \loe
u(x_\ell,y_\ell)$ or $u(x_\ell,y_\ell) \loe u(x_k,y_k)$. Since the intersection of $B(x_k,y_k)$ and
$B(x_\ell,y_\ell)$ is nonempty as they have the element $z$ in common, it follows that $B(x_k,y_k)\subs
B(x_\ell,y_\ell)$ or $B(x_\ell,y_\ell) \subs B(x_k,y_k)$. This proves that $(\cB(z),\subs)$ is narrow.

Take a ball $B(x,y)$ and $x',y'\in B(x,y)$, and set $\gamma:=u(x,y)$. Then $u(x,x')\leq \gamma$ and
$u(x,y')\leq \gamma$, hence by (U2), $u(x',y')\leq \gamma$. Assume that $u(x',y')=\gamma$, and take any $a\in
B(x,y)$, i.e., $u(x,a)\leq\gamma$. The latter together with $u(x,x')\leq \gamma$ implies that
$u(x',a)\leq \gamma=u(x',y')$, that is, $a\in B(x',y')$. Consequently, $u(x',y')=\gamma$ implies
$B(x,y)\subseteq B(x',y')$. Therefore, $B(x',y')\subsetneq B(x,y)\Rightarrow u(x',y')<u(x,y)$. Hence if
$\Gamma$ is countable, then $(\cB,\subs)$, and thus also $(\cB(z),\subs)$, admits only countable strictly
increasing sequences.
\end{pf}

\pars
Take an ultrametric space $(X,u)$ with countable narrow value set. Then the previous lemma shows that the
assumptions of Theorem~\ref{ThmHevyNarow} are satisfied, which proves Theorem~\ref{MT2}.

\section{A result on chain union stability}
\label{Sectnr}

Given a poset $\Gam$, define $\sig \Gam$ to be the set of all nonempty directed initial segments of $\Gam$ of countable cofinality. In other words, $C \in \sig \Gam$ if there exists a sequence $p_0 \loe p_1 \loe \cdots$ in $\Gam$ such that $C = \setof{x \in \Gam}{(\exists\, n) \; x \loe p_n}$. Then $\sig \Gam$ is a poset with inclusion that can be viewed as the ``$\omega$-completion'' of $\Gam$, namely, $\Gam$ is naturally embedded into $(\sig \Gam, \subs)$ and every countable chain in $\sig \Gam$ has a least upper bound (the union of that chain is directed and has a countable cofinality, by a standard diagonalization).

A poset $\Gam$ is \emph{up-countable} if for every $p \in \Gam$ the set $\setof{x \in \Gam}{p < x}$ is countable.

%
%

Given a ball space $(X, \cB)$, an \emph{ultradiameter} is a function $\map \delta \cB \Gam$, where $\Gam$ is a poset, such that $\delta$ is increasing with respect to inclusion and the following condition holds:
\begin{equation}
	(\forall \; B_0, B_1 \in \cB) \;\;\delta(B_0) \loe \delta(B_1) \Land B_0 \cap B_1 \nnempty \implies B_0 \subs B_1. \tag{{U}}
\end{equation}

\begin{tw}\label{THMboergoriehgi}
	Assume $(X,\cB)$ is a chain regular ball space with an ultradiameter $\map \delta \cB \Gam$ such that $\Gam$ is up-countable. Then $(X, \cB)$ is chain union stable.
\end{tw}

\begin{pf}
	Fix a chain $\Ef = \sett{F_\al}{\al<\kappa}$ in $\cu(\cB)$ such that $\al<\beta \implies F_\al \subsetneq F_\beta$, where $\kappa$ is an infinite cardinal. For each $\al<\kappa$ choose a chain $\Cee_\al \subs \cB$ with $\bigcup \Cee_\al = F_\al$.
	Choose $C_0 \in \Cee_0$. We may assume $C_0$ is minimal in $\Cee_0$.
	
	Using the chain regularity of $\cB$ and refining each of the chains $\Cee_\al$, we may assume that $C_0 \subs C$ for every $C \in \Cee_\al$, for every $\al<\kappa$.

	Let $\cB(C_0) = \setof{B \in \cB}{C_0 \subs B}$. So $\Cee_\al \subs \cB(C_0)$ for every $\al<\kappa$.
	Note that $\delta$ restricted to $\cB(C_0)$ satisfies
	$$\delta(B_0) \loe \delta(B_1) \iff B_0 \subs B_1.$$
	This is thanks to property (U).
	Hence $\delta \rest \cB(C_0)$ is an isomorphic embedding of $(\cB(C_0), \subs)$ into $\Gam$.
	Since $\Gam$ is up-countable, we conclude that $\cB(C_0)$ is countable. Lemma~\ref{LMerighworug} implies that $\kappa = \omega$. Finally, Lemma~\ref{LMfisgbuisgh} finishes the proof.	
\end{pf}

Note that more generally, the assertions of Lemma~\ref{LMetighueri} also hold when a ball space $(X,\cB)$ admits an
ultradiameter with values in a countable narrow poset, which is shown in the proof of
\cite[Theorem~3.4]{KKub}. 
Thus we obtain one more criterion for chain union stability.

\begin{tw}\label{THMPiecDwa}
	Assume $(X,\cB)$ is a ball space admitting an ultradiameter with values in a countable narrow poset. Then $(X,\cB)$ is chain union stable.
\end{tw}

\section{Final remarks and questions}


\begin{pyt}
	Can ``narrow'' be omitted from Thms~\ref{MT2}, \ref{ThmHevyNarow} and~\ref{THMPiecDwa}?
\end{pyt}

\begin{pyt}
	Can ``countable'' be omitted from Thm~\ref{THMPiecDwa}?
\end{pyt}

\end{document}